*To my daughters Vianne and Natalie*

# THE TRIANGLE ALTITUDES THEOREM IN HYPERBOLIC PLANE GEOMETRY

## NICHOLAS PHAT NGUYEN

**Abstract**. We provide a new formulation and proof of the triangle altitudes theorem in hyperbolic plane geometry, together with an easily computed discriminant to distinguish between different basic configurations of the altitudes of such a triangle.[1]

1. **INTRODUCTION**. In 2005, the late Vladimir Arnold discovered an algebraic proof of the triangle altitudes theorem in hyperbolic plane geometry using the Jacobi identity for Lie algebras. See [1] and [2]. In this paper, we present a new algebraic proof of the hyperbolic triangle altitudes theorem, using ideas from bilinear algebra. In addition, we also provide a readily computable discriminant that distinguishes all different altitude configurations for any hyperbolic triangle.

By an altitude, we mean a line passing through a vertex of the triangle and orthogonal to the opposite side. In Euclidean plane geometry, it is well-known that the three altitude lines of a triangle are concurrent. There is no such simple statement for a hyperbolic triangle. Instead, we have three different situations that together may be considered as the theorem about triangle altitudes in hyperbolic plane geometry. Given any triangle in a hyperbolic plane, we have the following three alternatives for the altitudes.

---

- Two of the altitudes intersect in a point in the hyperbolic plane, in that case the remaining altitude also passes through the same point.
- Two of the altitudes are asymptotically parallel (i.e., have a common point on the boundary at infinity), in that case the remaining altitude also is asymptotically parallel (i.e., passing through the same common point at infinity).
- Two of the altitudes are divergently parallel, i.e., they are orthogonal to a common line, and in that case the third altitude is also orthogonal to that common line. [2]

As another perspective, we will prove the following simple version of the triangle altitudes theorem in hyperbolic plane geometry.

**Theorem**: *The three altitudes of any hyperbolic triangle are linearly dependent.*

What do we mean by linear dependence? In the standard models of a hyperbolic plane, either based on the upper half plane or the interior unit disc in $\mathbf{R}^2$, the lines of that hyperbolic plane are given by Euclidean lines or circles in $\mathbf{R}^2$ that are orthogonal to the boundary horizontal line or boundary circle. Each Euclidean line or circle in $\mathbf{R}^2$ can be described as the zeros of a function $X \mapsto aX \cdot X + b \cdot X + c$, where $a$ and $c$ are real numbers, $b$ and $X$ are vectors in $\mathbf{R}^2$, and the dot product notation represents the standard scalar dot product in $\mathbf{R}^2$. The function or expression $aX \cdot X + b \cdot X + c$ is uniquely determined up to a scalar factor for each Euclidean line or circle, and the set $\mathcal{F}$ of all such functions or expressions has a natural structure as a real vector space of dimension 4, with component-wise addition. What our theorem says is that such an expression for any triangle altitude can be expressed as linear combination of the expressions for the other two altitudes, so that the three altitudes generate a subspace of dimension 2 of $\mathcal{F}$. (While the expression for a hyperbolic line is only determined up to a scalar factor, linear dependence or independence of a family of vectors is not changed by how we scale the individual vectors.)

---

[2] Some people also refer to the "divergently parallel" case as "ultra-parallel" and to the "asymptotically parallel" case as "parallel".

From that simple fact, we will deduce the different alternative configurations for triangle altitudes in the hyperbolic plane. In addition, we will provide a readily computable method to distinguish all alternative configurations for triangle altitudes. Specifically, given a hyperbolic triangle, we can calculate a number $\Delta$ that depends only on the expressions for the three sides of the triangle, such that:

- $\Delta > 0$ if and only if the three altitudes intersect in a common point.
- $\Delta = 0$ if and only if the three altitudes are asymptotically parallel.
- $\Delta < 0$ if and only if the three altitudes have a unique orthogonal line.

The prerequisite for this paper is relatively modest, requiring only a basic knowledge of symmetric bilinear forms, such as presented in the classic [3]. As for hyperbolic plane geometry, all the readers need to know is the basic definitions of hyperbolic lines in the upper half-plane model and what it means for two hyperbolic lines to be orthogonal, but we will recall those definitions at the appropriate time.

We add here a note on terminology. In general, a vector space of dimension 2 endowed with a nondegenerate symmetric bilinear form (a quadratic form) is either anisotropic or isometric to a standard space known as an Artinian or hyperbolic plane which has two linearly independent lines of isotropic vectors. Relative to a basis consisting the two linearly independent isotropic vectors, the matrix of the bilinear form is a 2-by-2 symmetric matrix that has zeros in the diagonal and a non-zero number in the cross diagonal. Over **R**, the matrix of such a bilinear form can also be put in the diagonal form [1, −1] relative to a suitable basis. The literature more commonly refers to such a space as a hyperbolic plane, but in order to avoid confusion, we will refer to such a space as an Artinian plane, which is another name that is sometimes used for it. See [4] at paragraph 13.1.4.4. The term "hyperbolic plane" often used in the literature for this type of quadratic space is rather undesirable, because the term hyperbolic plane can of course also mean a hyperbolic manifold of dimension 2 or a space that satisfies the axioms of hyperbolic plane geometry.

2. **THE SPACE OF CYCLES AND THE CYCLE PAIRING PRODUCT.** Consider the set $\mathcal{F}$ of all functions $p$ from $\mathbf{R}^2$ to $\mathbf{R}$ of the form $p(X) = aX.X + b.X + c$, where $X$ and $b$ are vectors in $\mathbf{R}^2$ and $a$ and $c$ are real numbers. The dot product notation of course refers to the standard scalar product $(x, y).(z, t) = xz + yt$ on $\mathbf{R}^2$. The set $\mathcal{F}$ can naturally be endowed with the structure of a real vector space of dimension $2 + 2 = 4$, being parametrized by the vector $b$ and the numbers $a$ and $c$. We will refer to such a function $p$ (if it is not the constant zero function) as a 2-cycle, 1-cycle, or 0-cycle depending on whether the degree of $p$ is 2 (coefficient $a$ is nonzero), 1 ($a$ is zero but $b$ is nonzero, or 0 (both $a$ and $b$ are zero).

All Euclidean lines in $\mathbf{R}^2$ are the sets of zeros of 1-cycles and all Euclidean circles in $\mathbf{R}^2$ are the sets of zeros of 2-cyles. We have used the name cycle rather than line or circle to emphasize the fact that we are looking at functions or expressions rather than the sets of zero points in $\mathbf{R}^2$ of such functions or expressions. Two expressions can define the same line or circle if they are scalar multiples of each other. Moreover, some cycles do not have any zeros, e.g., the 2-cycle $X.X+1$ has no zero in $\mathbf{R}^2$ and similarly for any 0-cycle $c$.

Aside from the natural structure of a vector space over $\mathbf{R}$, we can also endow $\mathcal{F}$ with a symmetric bilinear form $<\_,\_>$ as follows. Given $p = aX.X + b.X + c$ and $p^* = a^*X.X + b^*.X + c^*$, we define $<p,p^*>$ as $b.b^* - 2ac^* - 2a^*c$.

This mapping is clearly symmetric and bilinear. Moreover, it is non-degenerate, because $\mathcal{F}$ with this pairing is isometric to the orthogonal sum of $\mathbf{R}^2$ (with the pairing $b.b^*$) and an Artinian plane (with the pairing $-2ac^* - 2a^*c$). We will refer to this fundamental symmetric bilinear form on $\mathcal{F}$ as the cycle pairing or cycle product. By a suitable choice of basis, the matrix of the cycle pairing is diagonal of the form $[1, 1, 1, -1]$.

We can rewrite each 2-cycle $p$ as $aX.X + b.X + c = a(X + b/2a).(X + b/2a) - as$, where $s = (b.b - 4ac)/4a^2 = <p, p>/4a^2$. Note that $s$ has the same sign as $<p, p>$. By geometric analogy, we will refer to the point $(-b/2a)$ in $\mathbf{R}^2$ as the center of the 2-cycle $p$, regardless of whether or not the cycle has any zeros in $\mathbf{R}^2$. We will refer to $s$ as the square radius of the 2-cycle.

The following proposition summarizes the basic orthogonal relationships in $\mathcal{F}$ relative to the cycle pairing.

**Proposition 1**:

(a) Each 0-cycle is orthogonal to itself and all 1-cycles, but is not orthogonal to any 2-cycle. All 0-cycles are therefore isotropic.

(b) Two 1-cycles $p = \mathbf{b}.X + c$ and $q = \mathbf{d}.X + e$ have cycle product $<p, q> = \mathbf{b}.\mathbf{d}$. Consequently, $p$ and $q$ are orthogonal in $\mathcal{F}$ if and only if $\mathbf{b}$ and $\mathbf{d}$ are orthogonal in $R^2$ relative to the standard dot product. A 1-cycle is therefore non-isotropic.

(c) A 2-cycle $p$ with center $\mathbf{u}$ is orthogonal to a 1-cycle $q$ if and only if $q(\mathbf{u}) = 0$, i.e. the point $\mathbf{u}$ is a zero of the 1-cycle $q$, or that the Euclidean line defined by the 1-cycle $q$ passes through the point $\mathbf{u}$.

(d) Let $p = a(X + \mathbf{b}/2a).(X + \mathbf{b}/2a) - as$ and $q = d(X + \mathbf{e}/2d).(X + \mathbf{e}/2d) - dt$ both be 2-cycles. The cycle product of $p$ and $q$ is given by the following formula:

$$<p, q> = 2ad(s + t - w), \text{ where } w = (\mathbf{b}/2a - \mathbf{e}/2d).(\mathbf{b}/2a - \mathbf{e}/2d)$$

(e) A 2-cycle $p(X) = aX.X + \mathbf{b}.X + c$ is isotropic if and only if $\mathbf{b}.\mathbf{b} - 4ac = 0$, or equivalently, if $p(X) = a(X + \mathbf{b}/2a).(X + \mathbf{b}/2a)$. That is the expression for a Euclidean circle of zero radius centered at the point $(-\mathbf{b}/2a)$.

(f) Let $p = a(X + \mathbf{b}/2a).(X + \mathbf{b}/2a) - as$ and $q = d(X + \mathbf{e}/2d).(X + \mathbf{e}/2d) - dt$ be 2-cycles with $s \geq 0$ and $t \geq 0$. These expressions define Euclidean circles in $R^2$ with radius $\geq 0$. We have $<p, q> = 0$ if and only if the associated circles are orthogonal in the classical Euclidean sense.

(g) The only isotropic cycles in $\mathcal{F}$ are the 0-cycles and the 2-cycles corresponding to zero-radius Euclidean circles. The cycles $p$ of positive norm $<p, p>$ are the 1-cycles and the 2-cycles with positive square radius. These are the cycles that define Euclidean lines and circles with positive radius. The cycles $p$ of negative norm $<p, p>$ are the 2-cycles with negative square radius.

*Proof:*     Statements (a) and (b) are clear.  For statement (c), let $p = a(X + b/2a) \cdot (X + b/2a) - as = aX.X + b.X + c$ be a 2-cycle with center $u = -b/2a$, and let $q = e.X + f$ be a 1-cycle. We have $<p, q> = b.e - 2af = -2a(e.(-b/2a) + f) = -2a.q(u) = 0$ if and only if $q(u) = 0$. For statement (d), let:

$$p = aX.X + b.X + c = a(X + b/2a) \cdot (X + b/2a) - as, \text{ where } s = (b.b - 4ac)/4a^2$$

$$q = dX.X + e.X + f = d(X + e/2d) \cdot (X + e/2d) - dt, \text{ where } t = (e.e - 4df)/4d^2,$$

then by definition $<p, q> = b.e - 2af - 2dc$.

Because $w = (b/2a - e/2d) \cdot (b/2a - e/2d) = b.b/4a^2 + e.e/4d^2 - 2b.e/4ad$, we have, after cancelling out the terms $b.b/4a^2$ and $e.e/4d^2$:

$$s + t - w = -c/a - f/d + 2b.e/4ad, \text{ or}$$

$$2ad(s + t - w) = b.e - 2af - 2dc = <p, q>$$

For a 2-cycle $p$, $<p, p> = b.b - 4ac = 2a^2(2s - 0) = 0$ if and only if $s = 0$. That gives us statement (e).

For statement (f), note that $<p, q> = 0$ is equivalent to $s + t - w = 0$ in light of statement (d).  When $s$ and $t$ are both $\geq 0$,

$s =$ square of the radius of the circle associated with $p$
$t =$ square of the radius of the circle associated with $q$
$w =$ square of the distance between the centers of the two circles.

The equation $s + t - w = 0$ is then a simple restatement of the classical meaning of orthogonal circles in the Euclidean plane $\mathbf{R}^2$. In particular, a circle of zero radius is orthogonal to another circle of positive radius if and only if the center of the zero-radius circle lies on the other circle.

Statement (g) follows from the earlier statements.  ∎

Note that according to Proposition 1, Euclidean lines and circles are orthogonal in the classical Euclidean sense if and only if their corresponding cycles are orthogonal under the cycle pairing.

3.   **LINES AND ALTITUDES IN HYPERBOLIC PLANE**.   We will use the standard upper half plane model for our hyperbolic plane.  Specifically, our hyperbolic plane is the set $\mathcal{H}$ of all points $(x, y)$ in $\mathbf{R}^2$ with $y > 0$.

For each point $\boldsymbol{u}$ in the plane $\mathbf{R}^2$, we write $z(\boldsymbol{u}) = X.X - 2\boldsymbol{u}.X + \boldsymbol{u}.\boldsymbol{u}$ for the normalized 2-cycle centered at $\boldsymbol{u}$ with zero radius.  Note that $\boldsymbol{u}$ is a point or vector in $\mathbf{R}^2$, but $z(\boldsymbol{u})$ is a 2-cycle element of $\mathcal{F}$.  For distinct points $\boldsymbol{u}$ and $\boldsymbol{v}$ in the plane $\mathbf{R}^2$, we have $<z(\boldsymbol{u}), z(\boldsymbol{v})> = 4\boldsymbol{u}.\boldsymbol{v} - 2\boldsymbol{u}.\boldsymbol{u} - 2\boldsymbol{v}.\boldsymbol{v} = -2(\boldsymbol{u}-\boldsymbol{v}).(\boldsymbol{u}-\boldsymbol{v})$, which is always $< 0$.

From Proposition 1, note that a point $\boldsymbol{u}$ in the plane $\mathbf{R}^2$ lies on a Euclidean line or circle defined by a cycle $m$ if and only if $<z(\boldsymbol{u}), m> = 0$.

Consider the simple 1-cycle $p$ of the form $p(X) = \boldsymbol{j}.X$ where $\boldsymbol{j}$ is the vector $(0, -1)$ in $\mathbf{R}^2$.   The zero set of that 1-cycle is the horizontal line of all vectors $(x, y)$ whose second coordinate $y$ is zero.  Our hyperbolic plane $\mathcal{H}$ can be described as the subset of $\mathbf{R}^2$ consisting of all points $\boldsymbol{u}$ in $\mathbf{R}^2$ such that the product $<p, z(\boldsymbol{u})> = -2\boldsymbol{j}.\boldsymbol{u}$ is $> 0$.

The hyperbolic lines in the set $\mathcal{H}$ are defined as the portions in $\mathcal{H}$ of any Euclidean lines or circles in $\mathbf{R}^2$ that are orthogonal to the horizontal $x$-axis, i.e. the upper half of any vertical line and the upper semi-circle of any circle centered on the horizontal $x$-axis.

All these hyperbolic lines can be described as follows, in light of Proposition 1.

**Proposition 2**: *A hyperbolic line is a nonempty subset of $\mathcal{H}$ consisting of all points $\boldsymbol{u}$ in $\mathcal{H}$ such that $<m, z(\boldsymbol{u})> = 0$ for some nonzero cycle m orthogonal to p, i.e., $<m, p> = 0$.  If a cycle m defines a hyperbolic line, then the cycle m must have positive norm $<m, m>$.*

**Proof**.  When the hyperbolic line is a upper vertical ray, that hyperbolic line can be described as the set of points $\boldsymbol{u}$ in the upper half-plane that lie on the Euclidean line defined a 1-cycle $m$ (i.e., $<m, z(u)> = 0$) where $m$ is orthogonal to the horizontal $x$-axis represented by the cycle $p$.

When the hyperbolic line is an upper semicircle, that hyperbolic line can be described as the set of points $u$ in the upper half-plane that line on the Euclidean circle defined by a 2-cycle $m$ (i.e., $<m, z(u)> = 0$) of positive square radius centered on the horizontal $x$-axis (i.e., $<m, p> = 0$).

The above two cases cover all cycles that are orthogonal to both $p$ and a cycle $z(u)$ for some $u$ in $\mathcal{H}$. Indeed, we will show that such a cycle $m$ must have positive norm $<m, m>$, which means it must either be a 1-cycle or a 2-cycle with positive square radius.

We can easily eliminate the case of cycles with zero norm, because the only isotropic cycles are the zero-radius 2-cycles and the 0-cycles. If such a cycle is orthogonal to $p$, then it must be either $z(r)$ for some $r$ on the $x$-axis or a 0-cycle. None of these cycles can be orthogonal to some $z(u)$ with $u$ in the upper half plane.

Let us now see whether a cycle of negative norm can define a line in $\mathcal{H}$. Recall that with a suitable base of orthogonal vectors, the matrix of the cycle pairing has the diagonal form [1, 1, 1, –1]. A cycle of negative norm would give us a one-dimensional subspace of $\mathcal{F}$ isomorphic to the one-dimensional quadratic subspace [–1], whose orthogonal complement in $\mathcal{F}$ would be an anisotropic space of dimension 3 isomorphic to the positive definite quadratic space [1, 1, 1], according to the Witt's extension theorem. Because that complement is anisotropic, it contains no isotropic 2-cycle, and therefore there is no point $u$ such that $z(u)$ is orthogonal to the given cycle. So if a cycle $m$ has a zero point in $\mathcal{H}$, it must have positive norm. ∎

Note that if a cycle $m$ defines a line in $\mathcal{H}$, then any scalar multiple $\mu m$ of $m$ ($\mu \neq 0$) defines the same line. Recall that two hyperbolic lines are considered orthogonal, as a matter of definition, if and only if they are part of Euclidean lines or circles that are orthogonal in the classical Euclidean sense. In light of Proposition 1, that is equivalent to the cycles defining these hyperbolic lines being orthogonal under the cycle pairing.

**Proposition 3**: *Given a point **u** in $\mathcal{H}$ and a hyperbolic line L, there is one and only one hyperbolic line H through **u** and orthogonal to the line L.*

*Proof.* Let $\ell$ be a cycle that represents the hyperbolic line L. We note that the cycles $\ell$, p, and z(**u**) are linearly independent. Indeed, look at the 2-dimensional subspace of $\mathcal{F}$ generated by $\ell$ and p. With respect to that basis, the matrix of the cycle pairing on this subspace is diagonal (because $<\ell, p> = 0$), and because both $\ell$ and p have positive norms, the diagonal entries of the matrix are positive. That means on this subspace the cycle pairing is positive definite, and hence contains no isotropic cycle. Accordingly, the cycle z(**u**) is outside of that subspace, and therefore the subspace of $\mathcal{F}$ generated by $\ell$, p, and z(**u**) have dimension 3.

The orthogonal complement of that 3-dimensional subspace is therefore a one-dimensional subspace generated by some nonzero cycle h. By construction, this cycle is orthogonal to p, $\ell$, and z(**u**), and therefore defines a hyperbolic line H passing through **u** and orthogonal to the hyperbolic line L. This hyperbolic line H is unique because it must be defined by a cycle in the one-dimensional subspace generated by h, and proportional cycles define the same hyperbolic line. ∎

**Proposition 4**: *Suppose two different hyperbolic lines L and M intersect in a point **u** in $\mathcal{H}$. Let $\ell$ and m be two cycles that define the lines L and M. A hyperbolic line A passes through the same point **u** if and only if it can be defined by a nonzero linear combination of $\ell$ and m.*

*Proof.* Consider the 2-dimensional subspace U of $\mathcal{F}$ generated by p and z(**u**). This subspace is regular under the cycle pairing because $<p, z(\mathbf{u})> \neq 0$. Moreover, U is isotropic because it contains the isotropic cycle z(**u**). Accordingly, the cycle pairing on U is isometric to the Artinian plane [1, –1].

The orthogonal complement $U^\perp$ of U is a 2-dimensional subspace of $\mathcal{F}$ which contains $\ell$ and m. Because $\ell$ and m are linearly independent (because they define different hyperbolic lines), they form a basis of the subspace $U^\perp$. Now note that any hyperbolic line passing through **u** must be represented by a cycle in $U^\perp$, because a defining cycle for such a line must be orthogonal to both $p$ and $z(u)$. So if the line $A$ passes through $u$, it can be defined by a linear combination of $\ell$ and $m$.

Conversely, we will show that any nonzero linear combination of $\ell$ and m, i.e., any nonzero cycle in $U^\perp$, defines a hyperbolic line that passes through $u$. By the Witt's extension theorem, the space $U^\perp$ must be isometric to the positive definite plane [1, 1], and therefore any nonzero cycle in $U^\perp$ has positive norm. By Proposition 2, those cycles define hyperbolic lines if they are orthogonal to $p$. Moreover, those cycles are orthogonal to $z(u)$ by construction, and so they represent hyperbolic lines that pass through the point $u$. ■

4.  **TRIANGLE ALTITUDES IN HYPERBOLIC PLANE.** Consider a hyperbolic triangle defined by 3 points $u, v, w$ in $\mathcal{H}$. Let $L$ be the sideline of the triangle opposite $u$ (the line passing through $v$ and $w$), and similarly let $M$ and $N$ be the sidelines opposite $v$ and $w$ respectively. Let the lines $L, M,$ and $N$ be defined by cycles $\ell, m,$ and $n$ (which are determined only up to a scalar factor) respectively.

By Proposition 3, we know that through each vertex of the triangle there is a unique line that is orthogonal to the opposite sideline. Such a line is called an altitude line, or altitude for short.

Consider the altitude of our triangle through the vertex $u$. By Proposition 4, we know that such an altitude can be defined by a linear combination of $m$ and $n$ (because the lines $M$ and $N$ pass through $u$). Let $\alpha m + \beta n$ be a linear combination that defines this altitude through $u$. Because this altitude is orthogonal to the line $L$, we must have

$< \alpha m + \beta n, \ell > = \alpha <m, \ell> + \beta <n, \ell> = 0.$

That means we can choose $\alpha = <n, \ell>$ and $\beta = - <m, \ell>$.

So the cycle $A(u) = <n, \ell>m - <m, \ell>n$ represents the altitude through $u$. Similarly,

$A(v) = <\ell, m>n - <n, m>\ell$ represents the altitude through $v$, and

$$A(w) = \langle m, n \rangle \ell - \langle \ell, n \rangle m \quad \text{represents the altitude through } w.$$

We have $A(u) + A(v) + A(w) = 0$. Therefore the three altitudes are linearly independent and generate a 2-dimensional subspace $\mathcal{A}$ of $\mathcal{F}$.

Note that while each of the cycles $A(u)$, $A(v)$, and $A(w)$ is only determined up to a scalar factor, their linear dependence is a fact independent of any scaling of the individual cycles. We have therefore proved the

**Triangle Altitudes Theorem**: *The cycles representing the three altitudes of any hyperbolic triangle are linearly dependent.* ∎

Now consider the 2-dimensional subspace $\mathcal{A}$ of $\mathcal{F}$ that is generated by the three altitudes of the given triangle. On such a 2-dimensional subspace, we have the following disjoint alternatives for the cycle pairing.

*Alternative 1*: *The cycle pairing on $\mathcal{A}$ is anisotropic.*

We claim that this is the case if and only if two of the altitudes intersect in a point in $\mathcal{H}$. (In that case, the remaining altitude also passes through that point in light of our theorem and Proposition 4.)

Note that $\mathcal{A}$ is anisotropic if and only if the cycle pairing on $\mathcal{A}$ is positive definite (because the cycles in $\mathcal{A}$ have positive norm), which means the cycle pairing on $\mathcal{A}$ has diagonal form $[1, 1]$ relative to a suitable basis. But the subspace $\mathcal{A}$ is isometric to the quadratic space $[1, 1]$ if and only if its orthogonal complement in $\mathcal{F}$ is isometric to the Artinian plane $[1, -1]$ by Witt's extension theorem.

So our claim is equivalent to saying that two of the altitudes intersect at a point in $\mathcal{H}$ if and only if the orthogonal complement $\mathcal{A}^\perp$ of $\mathcal{A}$ is an Artinian plane.

If two of the altitudes, say the altitudes defined by $A(u)$ and $A(v)$, have an intersection point $s$, then $z(s)$ is a cycle orthogonal to both $A(u)$ and $A(v)$, and hence belongs to $\mathcal{A}^\perp$. $\mathcal{A}^\perp$

has dimension 2 and also includes $p$. $\mathcal{A}^\perp$ is therefore generated by $p$ and $z(s)$. It is regular (because $<p, z(s)> \neq 0$) and isotropic, and hence is an Artinian plane.

Now consider the converse. We want to show that if $\mathcal{A}^\perp$ is an Artinian plane, then all the lines defined by cycles in of $\mathcal{A}$ must intersect at a point in $\mathcal{H}$. Because $\mathcal{A}^\perp$ is isometric to an Artinian plane, it must contain two isotropic lines not orthogonal to $p$. These isotropic lines therefore must be multiples of the 2-cycles $z(s)$ and $z(t)$ for two points $s$ and $t$ in $\mathbf{R}^2$ that are off the horizontal $x$-axis. We want to show that one of the points $s$ or $t$ must be in the upper half-plane $\mathcal{H}$.

Consider the cycle $z(t) - 2 <p, z(t)>p$. It is a 2-cycle of leading coefficient 1, belongs to the same subspace as $z(t)$ and $p$, and is not a scalar multiple of $z(t)$. Moreover, that cycle is isotropic because its norm is equal to

$$- 2 <p, z(t)> <p, z(t)> - 2 <p, z(t)> <p, z(t)> + 4<p, z(t)> <p, z(t)> = 0$$

The cycle $z(t) - 2 <p, z(t)>p$ therefore must be $z(s)$. We therefore have

$$<p, z(s)> = <p, z(t) - 2 <p, z(t)>p> = <p, z(t)> - 2 <p, z(t)> = - <p, z(t)>.$$

It follows that one (and only one) of the points $s$ and $t$ are in $\mathcal{H}$. Let's say $s$ is in $\mathcal{H}$. By construction, $p$ and $z(s)$ are orthogonal to all the cycles in $\mathcal{A}$, meaning that all the nonzero cycles in $\mathcal{A}$ define lines that pass through $s$, including the three altitudes of the given triangle. Our claim is now proved.

Note that the cycle pairing on $\mathcal{A}$ is anisotropic if and only if the determinant of the cycle pairing matrix relative to any basis is positive.

*Alternative 2*: The cycle pairing on $\mathcal{A}$ is regular but isotropic.

We claim that this is the case if and only if there is a unique line in $\mathcal{H}$ orthogonal to all the lines defined by cycles in $\mathcal{A}$, including of course the three triangle altitudes.

The subspace $\mathcal{A}$ is regular and isotropic if and only if the cycle pairing on $\mathcal{A}$ has diagonal form [1, –1] relative to a suitable basis, which is the case if and only if the orthogonal complement $\mathcal{A}^\perp$ of $\mathcal{A}$ in $\mathcal{F}$ is isometric to the positive definite plane [1, 1], which is the case if and only if we can find a cycle $q$ in $\mathcal{A}^\perp$ that has norm 1 and is orthogonal to $p$ (which is a cycle of norm 1 in $\mathcal{A}^\perp$). Such a cycle $q$ is uniquely determined and therefore defines a unique hyperbolic line. By construction, that hyperbolic line is orthogonal to all the lines defined by cycles in $\mathcal{A}$, including the altitudes in question.

Note that the cycle pairing on $\mathcal{A}$ is regular and isotropic if and only if the determinant of the cycle pairing matrix relative to any basis is negative.

*Alternative 3*: *The cycle pairing on $\mathcal{A}$ is degenerate.*

We can eliminate the case when the cycle pairing on $\mathcal{A}$ is zero (i.e., the space $\mathcal{A}$ is totally isotropic). That is because the list of all isotropic cycles listed in Proposition 1 readily implies that $\mathcal{F}$ has no totally isotropic subspace of dimension 2. So what we are looking at is the case when the cycle pairing on $\mathcal{A}$ has a one-dimensional radical.

Note that the cycle pairing on $\mathcal{A}$ is degenerate if and only if the cycle pairing on its orthogonal complement $\mathcal{A}^\perp$ (which also has dimension 2) is degenerate, and that is the case if and only if $\mathcal{A} \cap \mathcal{A}^\perp$ is the common one-dimensional radical of $\mathcal{A}$ and $\mathcal{A}^\perp$. Let's say that this radical is generated by an isotropic cycle $r$. Because $<r, p> = 0$ (recall that the cycle $p$ is in $\mathcal{A}^\perp$), the fact that $r$ is an isotropic cycle means that it must be either: (i) a scalar multiple of $z(r)$ for some point $r$ on the horizontal $x$-axis; or (ii) a 0-cycle. We will regard all the 0-cycles as representing an exceptional point that we denote by the symbol ∞ and refer to as the extra point at infinity. The horizontal $x$-axis together with the point ∞ is regarded as the (completed) boundary at infinity of the hyperbolic plane $\mathcal{H}$.

Note that a 0-cycle is orthogonal to any 1-cycle, but never to a 2-cycle. So the point ∞ can be regarded as the common point at infinity of all the Euclidean lines defined by 1-cycles. Since the hyperbolic lines defined by 1-cycles are the upper vertical rays, all these vertical hyperbolic lines can be regarded as having a common point ∞ at infinity.

Based on the foregoing discussion, the cycle pairing on 𝓐 is degenerate if and only if all the cycles in 𝓐 are orthogonal to the isotropic cycle 𝓇, i.e., all the lines defined by cycles in 𝓐 pass through a common point on the boundary at infinity represented by the cycle 𝓇.

Note that the cycle pairing on 𝓐 is degenerate if and only if the determinant of the cycle pairing matrix relative to any basis is zero.

Therefore, we have three basic configurations for the triangle altitudes in a hyperbolic plane. To distinguish among these configurations, we can simply calculate the determinant Δ of the cycle pairing matrix on 𝓐 (relative to any basis) and see whether that determinant Δ is positive (alternative 1- anisotropic case), or negative (alternative 2 – non-degenerate but isotropic case), or zero (alternative 3 – degenerate case). Recall that if we change basis, then the determinant of the cycle pairing relative to the new basis is equal to the old determinant times a square scalar factor. Hence the sign of such a determinant is invariant.

Let's take *A(**u**)* and *A(**v**)* as a basis of 𝓐. Relative to this basis, we have

$$\Delta = \det \begin{pmatrix} <A(\mathbf{u}), A(\mathbf{u})> & <A(\mathbf{u}), A(\mathbf{v})> \\ <A(\mathbf{u}), A(\mathbf{v})> & <A(\mathbf{v}), A(\mathbf{v})> \end{pmatrix}$$

From the expressions for *A(**u**)* and *A(**v**)*, it's clear that Δ depends on just the cycles ℓ, 𝓂, and 𝓃 that define the 3 sides of the given triangle. If we adjust the cycles ℓ, 𝓂, and 𝓃 by scalar factors α, β, γ respectively, then *A(**u**)* and *A(**v**)* are each changed by the factor αβγ, and the determinant Δ is changed by the factor $(\alpha\beta\gamma)^4$, which of course does not change the sign of Δ.

Although our discussion has focused on the altitudes in a hyperbolic triangle, it applies more generally to any pair of distinct lines in the hyperbolic plane. Specifically, we have proved the following.

<u>**Configuration Theorem for Two Lines in Hyperbolic Plane**</u>: *Let L and M be two distinct lines in the hyperbolic plane 𝓗 defined by cycles ℓ and 𝓂. Let 𝓐 be the subspace of cycles* *generated by ℓ and m, and let* $\Delta = \det\begin{pmatrix} <\ell,\ell> & <\ell,m> \\ <\ell,m> & <m,m> \end{pmatrix}$ *be the determinant of the cycle pairing matrix on 𝓐 relative to the basis ℓ and m.*

*Then we have the following alternative situations:*

- *Δ > 0 if and only if the lines L and M intersect in a point in 𝓗.*
- *Δ = 0 if and only if the lines L and M are asymptotically parallel, i.e., if they have a common point on the boundary at infinity of 𝓗.*
- *Δ < 0 if and only if the lines L and M are divergently parallel, i.e., there is a unique line in 𝓗 orthogonal to both.*  ∎

NICHOLAS PHAT NGUYEN

*12015 12th Dr SE, Everett, WA 98208, U.S.A.*

*Email: nicholas.pn@gmail.com*